\documentclass{article}
\usepackage{amssymb}
\begin{document}

\bibliographystyle{plain}
\newcommand{\eqref}[1]{\mbox{\rm(\ref{#1})}}
\newcommand{\pl}{}
\newcommand{\C}[1]{{\protect\cal #1}}
\newcommand{\B}[1]{{\bf #1}}
\newcommand{\I}[1]{{\mathbb #1}}
\renewcommand{\O}[1]{\overline{#1}}
\newcommand{\binom}[2]{{#1\choose #2}}
\newcommand{\e}{\epsilon}
\newcommand{\qed}{\nolinebreak\mbox{\hspace{5 true pt}%
\rule[-0.85 true pt]{3.9 true pt}{8.1 true pt}}}

\newtheorem{theorem}{Theorem}
\newtheorem{lemma}[theorem]{Lemma}
\newtheorem{corollary}[theorem]{Corollary}
\newtheorem{proposition}[theorem]{Proposition}
\newtheorem{conjecture}[theorem]{Conjecture}
\newtheorem{example}[theorem]{Example}
\newtheorem{problem}[theorem]{Problem}

\newcommand{\MComp}{MComp}
\newcommand{\tp}{\tilde p}
\newcommand{\td}{\tilde d}
\newcommand{\xx}{{X\times X}}

\author{Oleg Pikhurko\thanks{Current address: DPMMS, Cambridge
University, Cambridge CB3 0WB, England. \mbox{E-mail:} {\tt
O.Pikhurko@dpmms.cam.ac.uk}}\\
 Department of Mechanics and Mathematics\\
 L'viv State Ivan Franko University\\
 L'viv 290602, Ukraine}

 \title{Operators Extending (Pseudo-)Metrics}
 \date{\small This paper appeared in the {\it Proccedings of All-Ukrainian
Conference of Young Researchers (Mathematics)}, {\bf N 1302}, Part~I (1994) 
10--16.\\
 This version (minor corrections): December 12, 2000}
 \maketitle

\begin{abstract}
 We introduce a general method of extending (pseudo-)metrics from
$X$ to $FX$, where $F$ is a normal functor on the category of
metrizable compacta. For many concrete instances of $F$, our method
specializes to the known constructions.\end{abstract}

\section{Introduction}

Consider the category of all compact metrizable spaces which will be
referred to as $\MComp$. All functors are expected to be {\em normal}
(for the definition and properties see~\cite[page 165]{2} or~\cite{3})
and to have $\MComp$ as both the domain and the codomain. For a normal
functor $F$, every space $X$ is naturally embeddable in $FX$, so further
in this work $X$ is considered to be a subspace of $FX$.

By an {\em operator} $u:C({-})\to C(F({-}))$ we mean a family of
maps
 $$(u_X:C(X)\to C(FX))_{X\in\MComp},$$
 where $C(X)$ denotes the set of all continuous mappings from $X$ to
$\I R$. Considering different topologies on this set, one can speak
about {\em operators continuous in the pointwise topology, in the
uniform topology,} etc. An operator is called a {\em functorial
operator} if for every $i:Y\to X$ the following identity holds:
 \begin{equation}\label{eq:1} u_Y\circ i_*=(F(i))_*\circ u_X.
\end{equation}
 Here, for $i:Y\to X$, the mapping $i_*:C(X)\to C(Y)$
corresponds $\phi$ to $\phi\circ i$. 

For $f,g\in C(X)$ we write $f\ge g$ to denote the poinwise inequality:
$f(x)\ge g(x)$ for all $x\in X$. An operator $u$ is an {\em extension
operator} if $u_X(\phi)|_X=\phi$; {\em monotonous} if $\phi\ge\psi$
implies $u_X(\phi)\ge u_X(\psi)$; {\em semiadditive} if
$u_X(\phi+\psi)\le u_X(\phi)+u_X(\psi)$; {\em positive} if $\phi\ge 0$
implies $u_X(\phi)\ge 0$, all $X$ and $\phi,\psi\in C(X)$.

Here we investigate a general method for extending
(pseudo-)metrics from a metrizable compact $X$ to $FX$, where $F$ is a
normal functor. For many concrete instances of $F$, our method
specializes to the known constructions.

\section{Definition and properties of the new operator}

Suppose that we have a normal functor $F$ and an operator $u:C({-})\to
C(F({-}))$. For $a,b\in FX$, $\langle a,b\rangle$ denotes the set
$$
\{c\in F(X\times X)\mid Fpr_1(c)=a,\
Fpr_2(c)=b\}=(F(pr_1),F(pr_2))^{-1}(a,b).
$$
It is not empty since any normal functor is bicommutative. Also, we
will use some other notation: 
 \begin{eqnarray}
&\Delta_X:X\to X\times X,\quad &\Delta(x)=(x,x);\\
&\nabla_X:X\times X \to X\times X,\quad &\nabla_X(x,y)=(y,x)
 \end{eqnarray}
If no confusion arises, we simply write $\Delta$ or $\nabla$.

For any real-valued function $p$ on $X^2$, we may define a function $\tp$
on $(FX)^2$ by the following formula:
 \begin{equation}\label{eq:2}
\tp(a,b)=\inf\{u_{X\times X}(p)(c)\mid c\in \langle a,b\rangle \},\quad a,b\in FX
 \end{equation}

The formula~(\ref{eq:2}) gives the promised operator~$\tilde{}$. Of
course, to define it, one needs an operator $u$ first, so it seems
that do not gain much. But, for many functors $F$, there is usually a
natural and obvious definition of $u$, while it is typically not clear how
to define a (pseudo-)metric on $FX$ should we have one on $X$.

\begin{lemma}\label{le:1}
If $u$ is an extension operator, then the function $\tp$ extends $p$.
\end{lemma}
 \smallskip{\it Proof.}  The claim is obvious because, for any normal functor $F$ and
arbitrary $a,b\in X$, the set $\langle a,b\rangle\subset F(X\times
X)$ consists of one point.\qed \medskip

\begin{lemma}\label{le:2}
If $u$ is a positive, monotonous, semiadditive functorial operator, then
for any pseudometric $p$ on $X$ the function $\tp$ is a pseudometric on
$FX$. 
\end{lemma}
 \smallskip{\it Proof.}  For any pair $(X\supset Y)$ and for every $\phi\in C(X)$ such that
$\phi|_Y=0$, we have $u_X(\phi)|_{FY}=0$. This can be deduced
from~\eqref{eq:1} by letting $i$ be the identity map $Y\to X$. Here,
$$ u_X(\phi)|_{FY}=(Fi)_*(u_X(\phi))=u_Y(i_*(\phi))=u_Y(0)=0.  $$

Now we can prove that, for any $a\in FX$, we have $\tp(a,a)=0$. Since
$p|_{\Delta(X)}=0$ we have $u_{X\times X}(p)|_{F\Delta(X)}=0$, and
$$
0\le \tp(a,a)\le u_{X\times X}(p)(F\Delta(a))=0.
$$

The function $\tp$ is symmetric, as
$$
 F\nabla(\langle b,a\rangle )=\langle a,b\rangle , \quad \forall a,b\in FX
$$
and
\begin{eqnarray*}
\tp(a,b)&=&\inf(u_{X\times X}(p)(\langle a,b\rangle))\ =\ \inf(u_{X\times
X}(p)(F\nabla(\langle b,a\rangle ))\\
 &=& \inf(u_{X\times X}(\nabla_*(p))(\langle b,a\rangle ))\ =\
\inf(u_{X\times X}(p)(\langle b,a\rangle )\ =\ \tp(b,a).
\end{eqnarray*}
 In this chain of equalities we used the symmetry of $p$ (i.e.\
$\nabla_*(p)=p$), the functoriality of $u$ (i.e.\ $u_{X\times
X}(\nabla_*(p))(x)=F\nabla_*\circ u_{X\times X}(p)(x)=u_{X\times
X}(p)(F\nabla(x)))$ and the identity
 \begin{equation}\label{eq:3}
F\nabla_*\circ u_{X\times X}=u_{X\times X}\circ\nabla_*:C(X\times X)\to
C(F(X\times X)).
 \end{equation}

Let $a$, $b$ and $c$ be arbitrary points in $FX$. Choose $x_1\in \langle a,b\rangle $
and $x_2\in \langle b,c\rangle $, such that $\td(a,b)=u_{X\times X}(x_1)$ and 
$\td(b,c)=u_{X\times X}(d)(x_2)$. $F$ is bicommutative so there exists
$y\in F(X^3)$ such that $Fpr_{12}(y)=x_1$ and $Fpr_{23}(y)=x_2$. Let
$x_3=Fpr_{13}(y)\in \langle a,c\rangle $. Then
\begin{eqnarray*}
\td(a,c)&\le& u_\xx(d)(x_3)\ =\
u_\xx(d)(Fpr_{13}(y))\\
 &=&
u_{X^3}(d\circ pr_{13})(y)\le u_{X^3}(d\circ pr_{12}+d\circ
pr_{23})(y)\\
 &\le& u_\xx(Fpr_{12}(y))+u_\xx(Fpr_{23}(y))\ =\ \td(a,b)+\td(b,c).
\end{eqnarray*}
 The lemma is proved.\qed \medskip

\begin{lemma}\label{le:3} If $u$ is continuous in the uniform
topology, then so is the operator~$\tilde{\ }$.  \end{lemma}
 \smallskip{\it Proof.}  For any $a,b\in FX$, we have
 \begin{eqnarray*}
 \|u_\xx(d_1)-u_\xx(d_2)\|_\infty &\ge&
u_\xx(d_1)(x_2,y_2)-u_\xx(d_2)(x_2,y_2)\\
 \ge\ \td_1(a,b)-\td_2(a,b)&\ge&
u_\xx(d_1)(x_1,y_1)-u_\xx(d_2)(x_1,y_1)\\
 &\ge& -\|u_\xx(d_1)-u_\xx(d_2)\|_\infty, \end{eqnarray*}
 where $\td_i(a,b)=u_\xx(d_i)(x_i,y_i)$, $i=1,2$. Hence
 $$ \|\td_1-\td_2\|_\infty\le \|u_\xx(d_1)-u_\xx(d_2)\|_\infty $$
 and the operator~$\tilde{\ }$ is continuous in the uniform topology.\qed \medskip

\begin{lemma}\label{le:4} If the mapping
 $$ H_X= (Fpr_1,Fpr_2):F(\xx)\to FX\times FX $$
 is open for any $X\in\MComp$, then $\td:FX\times FX\to \I R$ is
continuous.  \end{lemma}
 \smallskip{\it Proof.}  In fact, $\langle a,b\rangle =H_X^{-1}(a,b)$. Mapping $H_X$ is both
open and closed as $dom(H_X), codom(H_X)\in\MComp$. So the mapping
 $$ H_X^{-1}:FX\times FX\to \exp(F(\xx)) $$
 is continuous. Also, for any fixed $f\in C(X)$. the infimum map
$\inf_f:\exp(X)\to\I R$, defined by $\inf_f(A)=\inf f(A)$, is continuous.

Putting this all together we obtain the required.\qed \medskip

The direct consequence of Lemmas~\ref{le:1}--\ref{le:4} is the following.

\begin{theorem}\label{th:1} If $u_X$ is a positive, monotone,
semiadditive functorial operator extending functions from $X$ to $FX$,
then the operator $\tilde{\ }$ defined by formula~\eqref{eq:2} extends
pseudometrics from $X$ to $FX$.  Moreover, if $u$ is continuous in
the uniform topology, then so is the operator~$\tilde{\ }$; if $H_X$ is an
open mapping for all $X\in\MComp$, then the pseudometric $\td$ is
continuous for every continuous pseudometric $d$.\qed \end{theorem}

A remarkable fact about the above defined opeartor $\tilde{\ }$
is that in many cases it coincides with the well-known constructions,
as we are going to demonstrate now.

\section{Case $F=\exp$}

Let $F=\exp$ (the functor of all closed subsets equipped with the
Vietoris topology, see~\cite[page 139]{2}.) We define $u:C({-})\to
C(\exp({-}))$ by the formula $u_X(\phi)(A)=\sup(\phi(A))$, $\phi\in
C(X)$, $A\in \exp(X)$. 

\begin{theorem}\label{th:2}
For every metric $d$ on $X$, we have $\td=d_H$ (Hausdorff metric).
\end{theorem}
\smallskip{\it Proof.}  Let $A,B\in\exp(X)$, 
 $$M=d_H(A,B)=\inf\{\e>0\mid A_\e\supset B, B_\e\supset A\},$$
 where, for example, $A_\e=\{x\in X\mid d(x,A)\le \e\}$.

Then either there is $b\in B$ with $d(b,A)=M$ or there is $a\in A$
with $d(a,B)=M$. Since $pr_1(C)=A$ and $pr_2(C)=B$ for every $C\in
\langle a,b\rangle $, we have $u_\xx(d)(C)\ge M$, which implies $\td\ge
d_H$.

On the other hand, define
$$
C=\{(a,b)\in A\times B\mid d(a,B)=d(a,b)\mbox{ or } d(A,b)=d(a,b)\}.
$$
It is easy to prove that $C\in \langle a,b\rangle $ and
$u_\xx(d)(C)=M$.
Thus, we obtain that $\td=d_H$.\qed \medskip

\section{Case $F=({-})^n$.}

To define an operator $u$ one has to assign a certain number, given a
real-valued function $\phi$ on $X$ and a sequence $x_1,\dots,x_n\in
X$. It may be done in many ways but the following definitions are most
interesting: \begin{eqnarray*}
 u_\xx(\phi)(x_1,\dots,x_n)&=&\left(\sum_{i=1}^n\phi(x_i)^p\right)^{1/p},\quad
p\ge 1;\\
 u_\xx(\phi)(x_1,\dots,x_n)&=&\max_{i}(\phi(x_i)).
\end{eqnarray*}

The easy verification shows that corresponding operators $\tilde{\ }$
have the following appearence:
\begin{eqnarray*}
 \td(x,y)&=&\left(\sum_{i=1}^n d(x_i,y_i)^p\right)^{1/p};\\
 \td(x,y)&=&\max_i(d(x_i,y_i)).
\end{eqnarray*}

\section{Case $F=P$}

Let $P$ denote the functor of probability measures, see~\cite{1}.
The topology on the space $PX$ can be defined by means of the metric
$$
\bar d(\mu,\nu)=\inf\{\eta(d)\mid \eta\in P(\xx), Ppr_1(\eta)=\mu,
Ppr_2(\eta)=\nu\},\quad \mu,\nu\in PX
$$

Letting $u_X(\phi)(\mu)=\mu(\phi)$, $\mu\in PX$, $\phi\in C(X)$, one can
see that the definitions of $\bar d$ and $\td$ coincide.

\section{Case of the free (free abelian) group functor}

On the contrary to our default assumptions, here we suppose that the
functor $G({-})$, the free group functor, is defined on the category
of metrizable compacta with selected point. (The selected point plays
the role of the identity in $GX$.) 

The topology on the space $GX$ may
be defined in different ways. Among them are the constructions of
Swierczkowski and Graev.  To find distance between ``words'' $A,B\in
GX$ one has to find all {\em proper representations} $A=\prod_{i=1}^n
(a_i)^{\e_i}$ and $B=\prod_{i=1}^m(b_i)^{\sigma_i}$, $a_i,b_i\in X$,
$\e_i,\sigma_i=\pm 1$, that is, representations which have the same
number of letters and degrees coinciding exactly: $n=m$ and
$\e_i=\sigma_i$ for $1\le i\le n$. Then
 $$ d_1(A,B)=\inf\left(\sum_{i=1}^n d(a_i,b_i)\right), $$
 where the infimum is taken for all proper
representations. This is Graev's construction. That of Swierczkowski
(let us denote it by $d_2$) is nearly the same except we calculate the sum
only for all {\em different} pairs $(a_i,b_i)$.  Obviously, $d_1\ge
d_2$.

It turned out that these metrics can also be represented in the
form~\eqref{eq:2} for suitable $u$. Indeed, for $\phi\in C(X)$ and for
$A=\prod_{i=1}^n(a_i)^{\e_i}\in FX$ (written in the reduced form), let 
 $$ u_X(\phi)(A)=\sum_i\phi(a_i),$$
 but in the first case we take sum for all $i=1,\dots,n$ and in the
second for all different $a_i$'s. The points of the set $\langle
A,B\rangle $ are in the bijective correspondence with the proper
representations, which sends $C=\prod_{i=1}^n(c_i)^{\e_i}$
(in the reduced form) to the representations $A=\prod_{i=1}^n
pr_1(c_i)^{\e_i}$ and $B=\prod_{i=1}^n pr_2(c_i)^{\e_i}$. Since
 $$
u_\xx(d)(C)=\sum_id(pr_1(c),pr_2(c))
$$
 we get the claimed result.

The case $F=A$ (the free abelian group functor) is analogous. The interested
reader should be able to transfer easily all results by himself.

\section{Acknowledgements}

I am very grateful to Michael Zarichnyi who suggested this
construction for my Diploma Work at the L'viv State University and who
kindly offered other help on this study.

\end{document}